\documentclass[11pt]{article}

\usepackage{amsbsy,amsfonts,amsmath,amssymb,eucal,mathrsfs}
\usepackage[all]{xy}

\def\no{\noindent}
\def\itemn#1{\item[\hspace{0.6mm} {\rm (#1)}]}

\def\et{{\rm et}}

\def\red{{\rm red}}
\def\tir{\mbox{-}}
\def\ev{{\rm ev}}
\def\tr{{\rm tr}}
\def\disc{{\rm disc}}

\renewcommand{\bar}{\overline}

\newtheorem{definition}{Definition}[subsection]

\newtheorem{prop}[definition]{Proposition}
\newtheorem{lemm}[definition]{Lemma}

\newtheorem{theo}[definition]{Theorem}
\newtheorem{construction}[definition]{Construction}

\newtheorem{remark}[definition]{Remark}

\newtheorem{remarks}[definition]{Remarks}

\newtheorem{question}[definition]{Question}
\newenvironment{ques}{\begin{question} \rm}{\end{question}}
\newtheorem{example}[definition]{Example}
\newenvironment{exam}{\begin{example} \rm}{\end{example}}
\newtheorem{examples}[definition]{Examples}

\newtheorem{nothing}[definition]{$\!\!$}
\newenvironment{noth}{\begin{nothing} \rm}{\end{nothing}}
\newenvironment{proo}{{\flushleft \bf Proof :}}{\hfill $\square$ \vspace{5mm}}
\newenvironment{prooU}{{\flushleft \bf Proof of unicity :}}{\hfill
  $\square$ \vspace{5mm}}
\newenvironment{prooE}{{\flushleft \bf Proof of existence :}}{\hfill
  $\square$ \vspace{5mm}}

\DeclareMathOperator{\Id}{id}
\DeclareMathOperator{\pr}{pr}
\DeclareMathOperator{\Spec}{Spec}
\DeclareMathOperator{\Sym}{Sym}

\DeclareMathOperator{\Hom}{Hom}
\DeclareMathOperator{\End}{End}
\DeclareMathOperator{\Alg}{Alg}

\DeclareMathOperator{\Frac}{Frac}
\DeclareMathOperator{\GL}{GL}
\DeclareMathOperator{\MM}{M}
\DeclareMathOperator{\HH}{H}
\DeclareMathOperator{\TT}{T}

\DeclareMathOperator{\CH}{CH}
\DeclareMathOperator{\Funct}{Funct}
\DeclareMathOperator{\Sch}{Sch}
\DeclareMathOperator{\Res}{Res}

\def\cA{{\cal A}}   
\def\cF{{\cal F}}   
   
\def\cO{{\cal O}}   \def\cS{{\cal S}}

\def\sA{{\mathscr A}} \def\sB{{\mathscr B}} \def\sC{{\mathscr C}}

\def\bA{{\mathbb A}}   \def\bF{{\mathbb F}}
\def\bG{{\mathbb G}}  
  \def\bZ{{\mathbb Z}}

 \def\fA{\mathfrak{A}} \def\fB{\mathfrak{B}}

\def\fAlg{\mathfrak{Alg}}

\begin{document}

\begin{center}
{\bf \Large Determinants of finite dimensional algebras}

\bigskip
\bigskip

{\bf Matthieu Romagny}
\end{center}

\bigskip

\small
\begin{center}
\begin{minipage}{15cm}
\noindent {\bf Abstract}: To each associative unitary
finite-dimensional algebra over a normal base, we associate a canonical
multiplicative function called its {\em determinant}. We give various
properties of this construction, as well as applications to the
topology of the moduli stack of $n$-dimensional algebras.
\end{minipage}
\end{center}
\normalsize

{\def\thefootnote{\relax}
\footnote{ \hspace{-6.8mm}
Date~: December 4, 2007. }}

\vspace{-5mm}

\section{Introduction}

Our object of interest in this article is the moduli stack of
$n$-dimensional associative algebras with unit, denoted
$\fAlg_n$. Given a free module with basis $e_1,\dots,e_n$, an
algebra structure is given by the constants $c^{ij}_k$
such that $e_ie_j=\sum_k\,c^{ij}_ke_k$, satisfying the relations
coming from the associativity rule
$(e_ie_j)e_k=e_i(e_je_k)$. Therefore, algebras with a given basis are
classified by the affine variety of structure constants
$c^{ij}_k$ cut out by these relations, denoted $\Alg_n$. There is an
action of the group $G$ of base changes, and $\fAlg_n$ is the quotient
stack of $\Alg_n$ by this action. The geometry of $\Alg_n$ remains rather
mysterious, and the attention of the specialists has focused
on topics like the determination of the number of irreducible components (open
question, first raised by Gabriel in \cite{Ga}), asymptotic bounds for
their dimensions, and the existence of smooth components (\cite{LBR})
or nonreduced components (\cite{DS}). For more
details and references, we refer to \cite{LBR} and \cite{P}.

In this paper, we introduce a chain of closed substacks
$$\fAlg_{n,\le 2}\subset \fAlg_{n,\le 3}\subset\dots\subset\fAlg_{n,\le
  n}=\fAlg_n
$$
indexed by the {\em degree of algebraicity $d$}. Roughly, an algebra
is in $\fAlg_{n,\le d}$ if all its monogenic subalgebras
have rank less than $d$. For example, by the Cayley-Hamilton theorem, the
algebra of $(n,n)$ matrices lies in $\fAlg_{n^2,\le n}$. Then, we construct
a {\em determinant} on the normalization of the reduced substack of
$\fAlg_{n,\le   d}\setminus \fAlg_{n,\le d-1}$, by which we mean a function
enjoying the properties listed in the result below which is our main theorem.

Before we state it, a word about the terminology is necessary. Any
locally free algebra over a ring~$R$ (or locally free sheaf of
algebras over $S=\Spec(R)$) may be seen as a vector bundle over $S$
endowed with an algebra structure over the ring scheme $\bG_{a,S}$ ({\it cf}
section~\ref{alg_and_alg_schemes}).
We use the point of view of algebra schemes because most identities
that we prove ({\it e.g.} multiplicativity of the determinant) are
polynomial identities and not just equalities of functions. Also this
choice makes our statements at the same time more precise and more
elegant. However, the reader unfamiliar with the vocabulary of schemes
may very well replace all ``algebra schemes over a scheme $S$'' by
``(ordinary) algebras over a ring $R$''. In the text, a special effort
is made to translate into commutative algebra the statements with a
strong flavour of algebraic geometry, wherever it seems necessary.
Now here is our main result~:

\bigskip

\no {\bf Theorem~(\ref{det}) } {\em Let $S$ be a normal integral
  scheme, $A/S$ an $n$-dimensional $\bG_{a,S}$-algebra scheme with
  $n\ge 2$, $A^\vee/S$ its linear dual, and~$d\ge 2$ the least integer such
  that $A/S$ belongs to $\fAlg_{n,\le d}$. Then there exists a unique
  section of $\Sym^d(A^\vee)$
denoted $\det_{A/S}$ or simply $\det$, such that~:
\begin{trivlist}
\itemn{1} $\det\colon A\to\bG_{a,S}$ is a morphism of multiplicative
unitary monoid schemes, and
\itemn{2} the morphism of schemes $A\to A$ defined by evaluating on
$a\in A$ the polynomial $P(T)=\det(T-a)$ is zero. (Here $T$ is a
scalar in the $\bG_{a,S}[T]$-algebra $A\times_S\bG_{a,S}[T]$.)
\end{trivlist}
The determinant satisfies the further properties~:
\begin{trivlist}
\itemn{3} the group scheme of units of $A$ is the
preimage of the multiplicative group $\bG_{m,S}$ under $\det$~;
\itemn{4} the formation of $\det$ is compatible with flat extensions of
normal integral schemes $S'\to S$.
\end{trivlist}}

\bigskip

The existence of the determinant has some consequences on the topology
of $\fAlg_n$. In fact we hope that it helps to sort out the irreducible
components (see subsection~\ref{appl_top} and in particular
question~\ref{ques_irred_comp}). As an example, the irreducible
component of $\MM_n(k)$ has received special attention due to its
connection with moduli spaces of vector bundles on curves
(see~\cite{LBR}). The following proposition, proved using the
determinant, implies that it is contained in $\fAlg_{n^2,\le n^2-n}$~:

\bigskip

\no {\bf Proposition~(\ref{an_application})}
Let $n_1,\dots,n_r\ge 1$ be integers and $n=(n_1)^2+\dots+(n_r)^2$. Assume
that one of the $n_i$ is at least $2$ and denote by $\nu$ their infimum.
Then the irreducible component of the algebra
$\sA_0=\MM_{n_1}(k)\times\dots\times\MM_{n_r}(k)$ is contained in
$\fAlg_{n,\le n-\nu}$.

\bigskip

Let us finally survey the contents of the paper. In
section~\ref{alg_and_alg_schemes}, we explain the correspondance between
algebras, or sheaves of algebras, and algebra schemes. We include the case of
infinite dimension to incorporate the polynomial algebra
$\bG_{a,S}[T]$ that appears in the main theorem. In section~\ref{sec_deg_alg}
we define the degree of algebraicity and the closed
substacks $\fAlg_{n,\le d}$, and we give their basic properties. In
section~\ref{determinants} we prove the main result (theorem~\ref{det}) and we
use it to define the determinant on the normalization of the reduced substack
of $\fAlg_{n,\le d}\setminus \fAlg_{n,\le d-1}$
(proposition~\ref{det_on_Alg_d}~; we conjecture that it is actually
unnecessary to normalize). We also prove that the determinant is
invariant under all (anti)automorphisms of the algebra
(proposition~\ref{invariance_of_det}). In section~\ref{comp_appl} we provide
basic computations of determinants~: determinant of the opposite algebra
(\ref{det_of_A_opp}), determinant of a product (\ref{product}). Also we
use the determinant to define intrinsic invariants of an algebra~: trace,
discriminant, unimodular group (\ref{tr_and_disc}) and to study the topology
of $\fAlg_n$ (\ref{appl_top}). Finally in section~\ref{examples}
we compute more examples of determinants~: field extensions,
quaternion algebras, exterior algebras, and three-dimensional algebras.


\tableofcontents

\section{Algebras and algebra schemes} \label{alg_and_alg_schemes}

In this paper, all algebras are unitary and associative, but not necessarily
commutative.

\subsection{Different ways to see an algebra}

\begin{noth}
Let $S$ be a scheme. If $\cF$ is a quasi-coherent $\cO_S$-module, we
denote the dual module by~$\cF^\vee$. Associated to $\cF$ are the {\em tensor
  algebra} $\TT(\cF)$, the {\em symmetric algebra} $\Sym(\cF)$, and the {\em
  exterior algebra} $\wedge\cF$. If $\cF=(\cO_S)^n$ for some integer $n$, and
we denote by $U_1,\dots,U_n$ its canonical basis, then $\TT(\cF)$ is the
noncommutative polynomial algebra $\cO_S\{U_1,\dots,U_n\}$ and $\Sym(\cF)$ is
the commutative polynomial algebra $\cO_S[U_1,\dots,U_n]$. More generally, if
$\cF$ is locally free of finite rank, then $\TT(\cF)$ and $\Sym(\cF)$ are
twisted (non)commutative polynomial algebras.

Let $\cA$ be a quasi-coherent $\cO_S$-algebra and $a=(a_1,\dots,a_n)$ a tuple
of sections of $\cA$ over~$S$. The morphism of {\em evaluation at $a$}
is the morphism of algebras
$$\ev_a\colon\cO_S\{U_1,\dots,U_n\}\to\cA$$
(or $\ev_a\colon\cO_S[U_1,\dots,U_n]\to\cA$ if the sections $a_i$ commute)
defined by $P(U_1,\dots,U_n)\mapsto P(a_1,\dots,a_n)$. Alternatively, in
defining these notions, one can exhibit a more canonical domain. Namely, one
can define $\cF$ to be the sub-$\cO_S$-module of $\cA$ generated by
$a_1,\dots,a_n$, and consider the morphism of algebras $\TT(\cF)\to\cA$,
resp. $\Sym(\cF)\to\cA$, induced by the inclusion $\cF\to\cA$. This map may again
be called {\em evaluation at $a$}~; note that $\TT(\cF)$, resp. $\Sym(\cF)$,
is indeed a polynomial algebra only if $\cF$ is free. Finally we say that {\em
  $a$ generates $\cA$} if $\ev_a$ is surjective.
\end{noth}

\begin{noth}
It is sometimes possible to view an $\cO_S$-module or $\cO_S$-algebra as a
scheme. In order to explain this, we need to recall that, just in the
same way as one defines the notion of a group scheme, we have the notion of a
ring scheme. The most important example of a commutative ring scheme is the
{\em structure ring scheme} denoted $\bG_{a,S}$, or simply $\bG_a$ if no
confusion seems likely to result. Its underlying scheme is the affine line
over $S$, and it has two composition laws $\bG_a\times_S\bG_a\to\bG_a$ called
addition and multiplication, and two sections $S\to \bG_a$ denoted $0,1$ which
are neutral for these laws.

In the same way as one identifies an $S$-scheme with its functor of
points, we will identify an $\cO_S$-algebra (or an $\cO_S$-module)
with its functor of sections. This is made possible by the embedding
$\cS ec\colon (\cO_S\tir\Alg)\to \Funct((\Sch/S)^\circ,\bG_{a,S}\tir\Alg)$
defined as follows~: the functor $\cS ec$ takes $\cA$ to
the functor $\cS ec_\cA$ defined
by
$$
\cS ec_\cA(T\stackrel{f}{\to}S)= \mbox{ the }
\HH^0(T,\cO_T)\mbox{-algebra } \HH^0(T,f^*\cA)\;.
$$
The same works for modules.

The functor of sections of an $\cO_S$-algebra $\cA$ whose underlying
$\cO_S$-module is locally free of finite rank $n$ is (representable by) a
$\bG_{a,S}$-algebra scheme whose underlying $\bG_{a,S}$-algebra module is an
$n$-dimensional vector bundle. Indeed, it
is representable by $A=\Spec(\Sym(\cA^\vee))$. Concretely, if $U=\Spec(R)$ is
an open affine subscheme of $S$ such that $\cA(U)$ is free, then given a basis
$e_1,\dots,e_n$ of $\cA(U)$, there are some constants $c^{i,j}_k\in R$ such
that $e_ie_j=\sum_k\,c^{i,j}_ke_k$.
The forms $t_i=e_i^*$ are a system of coordinates for~$A$, that
is to say $A=\Spec(R[t_1,\dots,t_n])$, and the multiplication of
$A$ is given, on the level of functions, by $t_k\mapsto
\sum_{i,j}\,c^{i,j}_k t_i\otimes t_j$. For simplicity, in all the paper we
will call such an algebra an {\em $n$-dimensional algebra}.

We take the opportunity here to mention that our conventions for vector
bundles differ slightly from the ones in~\cite{EGA}, in that we will
use the words
{\em vector bundle} only for those bundles whose sheaf of sections is locally
free. Also, if $\cF$ is a locally free module of finite rank on a scheme $S$,
then the {\em associated vector bundle} will be
$F=\Spec(\Sym(\cF^\vee))$ and not $\Spec(\Sym(\cF))$.

More generally, an $\cO_S$-algebra whose underlying module is locally
free of arbitrary rank is mapped to an inductive limit of
finite-dimensional vector bundles over $S$~; in the inductive system, the
morphisms are vector bundle homomorphisms.
Examples include the polynomial algebras $\bG_{a,S}\{U_1,\dots,U_n\}$ and
$\bG_{a,S}[U_1,\dots,U_n]$ associated to the polynomial algebra sheaves
$\cO_S\{U_1,\dots,U_n\}$ and $\cO_S[U_1,\dots,U_n]$.
\end{noth}

\subsection{Universal elements} \label{universal_elements}

In order to establish some identities concerning algebra schemes, it will
often be enough to check these identities for their
{\em universal element}. Let us indicate briefly the relevant notations and
terminology.

\begin{noth} \label{image_of_univ_elt}
Let $A\to S$ be a morphism of schemes. The {\em universal element} of $A/S$
is the section of the pullback scheme $\pr_2\colon A\times_S A\to A$ which is
given by the diagonal of $A/S$. It is denoted by the letter $\alpha$.
If $f\colon A\to B$ is a morphism of schemes over $S$, then
$(f\times\Id_A)\circ\alpha=f^*\beta$ where $\alpha$ and $\beta$ are the
universal elements of $A$ and $B$. If $f$ is scheme-theoretically dominant,
that is to say if the scheme-theoretic image of $f$ is $B$, then this means
that the universal element of $A$ maps to the universal element of $B$.
\end{noth}

\begin{noth}
For local computations, when $A$ is an $n$-dimensional $\bG_{a,S}$-algebra
scheme and $S$ is the spectrum of a ring $R$, it is useful to have a
description in terms of commutative algebra. In this case $A$ is determined by
its algebra of global sections $\sA:=A(R)$. We denote by
$R_\sA:=\Sym(\sA^\vee)$ its function ring, which is a graded $R$-algebra~;
here $\sA^\vee$ is the dual $R$-module of $\sA$.
Because $\sA$ is projective of finite rank, there is an isomorphism of
$R$-modules $\Hom_R(\sA,\sA)\simeq \sA\otimes_R \sA^\vee$, and since $R_\sA$
contains $\sA^\vee$ as the piece of degree $1$, the target module is a
submodule of $\sA\otimes_R R_\sA$.
In this setting, the universal element $\alpha$ is the image of the identity
$\Id\in\Hom_R(\sA,\sA)$ inside $\sA\otimes_R R_\sA$.
If $\sA$ is free as an $R$-module, and $e_1,\dots,e_n$ is a basis, then $R_\sA$
is the (commutative) polynomial $R$-algebra on the independant
variables $t_i=e_i^*$ and $\alpha$ is just $t_1e_1+\dots+t_ne_n$.

Now let $A,B,C$ be finite-dimensional $\bG_{a,S}$-algebra schemes over
$S=\Spec(R)$, with sets of global sections $\sA$, $\sB$, $\sC$ and
function rings $R_\sA$, $R_\sB$, $R_\sC$.

If $f:A\to B$ is a surjective homomorphism, there is an associated map
$\sA\to\sB$ and an injection
$R_\sB\hookrightarrow R_\sA$. It is easy to check that the image of the
universal element $\alpha\in A\otimes_RR_\sA$ under
the map $f\otimes\Id\colon A\otimes_R R_\sA\to B\otimes R_\sA$ is the
universal element $\beta\in B\otimes_R R_\sB\hookrightarrow B\otimes_R R_\sA$.
Indeed, this assertion is local on $R$, so we may localize and choose a basis
$e_1,\dots,e_n$ of $A$ such that $f$ is the projection onto the subspace
spanned by the first $k$ vectors $e_1,\dots,e_k$. Then
$\alpha=\sum_{i=1}^n\,e_i^*e_i$ and $f(\alpha)=\sum_{i=1}^k\,e_i^*e_i=\beta$.

In particular, a product algebra $C=A\times B$ has function ring
$R_\sC=R_\sA\otimes_R R_\sB$, with universal element $\gamma=(\gamma_1,\gamma_2)\in
\sC\otimes_R R_\sC$. By the above, the projection~$\gamma_1$ in $\sA\otimes_R R_\sC$ is
$\alpha\in \sA\otimes_R R_\sA\subset \sA\otimes_R R_\sC$, and the projection
$\gamma_2$ in $\sB\otimes_R R_\sC$ is $\beta$. We point out that if $A=B$,
it is not true that $\gamma_1=\gamma_2$, because $R_\sA=R_\sB$ is embedded in
two different ways in $R_\sC$.
\end{noth}

\section{Degree of algebraicity} \label{sec_deg_alg}

\subsection{Definition} \label{deg_alg}

Let $A\to S$ be an $n$-dimensional $\bG_{a,S}$-algebra scheme.
We say that {\em the degree of algebraicity of $A/S$ is less
  than or equal to~$d$}, written $\deg(A/S)\le d$, if the morphism~:
$$
\begin{array}{rcl}
A & \to & \wedge^{d+1}A \\
a & \mapsto & 1\wedge a\wedge\dots\wedge a^d \\
\end{array}
$$
is the zero morphism. This is the same thing as saying that the universal
element $\alpha$ satisfies $1\wedge \alpha\wedge\dots\wedge \alpha^d=0$.
Locally on $S$, we may choose a basis for $A$, and then the latter
condition is equivalent to the vanishing of the minors of size $d+1$
of the matrix whose columns are the powers of $\alpha$.
If $S$ is the spectrum of a field $k$, then we say that {\em the degree of
  algebraicity of $A/k$ is equal to~$d$}, written $\deg(A/k)=d$, if it
is less than $d$ but not less than $d-1$. If $S$ is arbitrary, we say 
that $\deg(A/S)=d$ if $\deg(A/S)\le d$ and all the fibres of
$A\to S$ have degree of algebraicity equal to~$d$.
If $A/S$ has a degree of algebraicity $d$, then $A\times_S S'$ has
degree of algebraicity $d$ for all extensions $S'\to S$.
Also if $A$ has a degree of algebraicity $d$, then $d\le n$ because
$\wedge^{n+1} A=0$. To better illustrate this definition, consider the
following example.

\begin{exam} \label{exam1}
Let $R=k[\epsilon]/\epsilon^2$ be the ring of dual numbers over a field $k$.
Consider the $3$-dimensional commutative algebra $\sA=R[x,y]/(x^2-\epsilon
x,xy,y^2)$ and its associated $R$-algebra scheme $A$. If we denote the
universal element by $\alpha=r+sx+ty$, we
have $1\wedge\alpha\wedge\alpha^2=-\epsilon s^2t\cdot 1\wedge x\wedge y$. This
is zero modulo~$\epsilon$, so $\deg(A\otimes k)=2$, but we do not attribute a
degree of algebraicity to $A$ itself.
\end{exam}

\subsection{Relation with monogenic subalgebras}

If $k$ is a field, then any element $x$ in a finite-dimensional
(ordinary) $k$-algebra has a minimum polynomial, whose degree is also
the dimension of the subalgebra generated by $x$~; we call it the {\em
  degree of $x$}, denoted $\deg(x)$. The degree of algebraicity of an
algebra as defined above has a simple meaning in these terms~:

\begin{lemm} \label{monogenic}
Let $A/S$ be an $n$-dimensional algebra scheme that has a degree of
algebraicity~$d$.
\begin{trivlist}
\itemn{1} If $S$ is the spectrum of a field $k$, there is a finite
separable field extension $k'/k$ such that $d=\sup\,\{\deg(x),\,x\in
A(k')\}$. If $k$ is infinite we may take $k'=k$.
\itemn{2} We have $d=n$ if and only if there exists a surjective {\'e}tale
extension $S'\to S$ such that $A\otimes_S S'$ is monogenic. If $S$ has
infinite residue fields, we may take for $S'$ a Zariski open cover
of $S$. In particular, if $d=n$ then $A$ is commutative.
\end{trivlist}
\end{lemm}

\begin{proo}
(1) Choose a basis $e_1,\dots,e_n$ for $\sA=A(k)$ and let $t_i=e_i^*$, so
$\alpha=t_1e_1+\dots+t_ne_n$. Also let $M_i$ denote the matrix with columns
  $1,\alpha,\dots,\alpha^i$. To say that $d=\deg(A)$ means that all the
minors of size $d+1$ of $M_d$ vanish, and one of the minors of size $d$ of
$M_{d-1}$ does not vanish, call $m\in k[t_1,\dots,t_n]$ this minor. Thus,
clearly, all elements of $\sA$ have degree less than $d$,
since they are specializations of $\alpha$. Now if $k$ is infinite, there is a
tuple of elements $(x_1,\dots,x_n)\in k^n$ where $m$ takes a nonzero value, so
$x=x_1e_1+\dots+x_ne_n$ has degree $\deg(x)=d$. If $k$ is finite there is such a
tuple in an algebraic closure of $k$, hence in a finite separable extension
$k'/k$.

\no (2) The {\em if} part is obvious, we focus on the {\em only
  if}. The claim is local on $S$ so we may assume that $S=\Spec(R)$ is
affine and $A$ is trivial as a vector bundle, determined by
$\sA=A(R)$. Let~$m$ be a maximal ideal of $R$ and $k=R/m$. By point~(1)
there is a finite separable field extension $k'/k$ such that
$\sA\otimes_R k'$ is monogenic, and if $k$ is
infinite we may take $k'=k$. Let $P\in R[T]$ be a monic polynomial with
reduction $p\in k[T]$ such that $k'=k[T]/(p)$, and let $R'=R[T]/(P)$.
Since $R'/R$ is finite flat and {\'e}tale at $m$, we may localize in $R'$ if
necessary and assume that $R'/R$ is {\'e}tale. Now let $g\in \sA\otimes_R R'$
be a lift of a generator of $\sA\otimes_k k'$.
By Nakayama's lemma, after a further localization in $R'$, $g$ is a generator
of $\sA\otimes_R R'$.
The last statement about commutativity follows immediately.
\end{proo}

\begin{exam}
Consider the commutative $\bF_2$-algebra
$\sA=\bF_2[x,y]/(x^2-x,xy,y^2-y)$ and its associated $\bF_2$-algebra scheme
$A$. Then $\sA$ satisfies $x^2=x$ for all $x\in \sA$, however its degree
is not $2$ because the morphism
$$
\begin{array}{rcl}
A & \to & \wedge^3A \\
a & \mapsto & 1\wedge a\wedge a^2 \\
\end{array}
$$
is $a=u+vx+wy\mapsto (vw^2-v^2w)(1\wedge x\wedge y)$ which is not zero.
In fact $\sA$ has degree $3$ and indeed $\sA\otimes \bF_4$ is generated as an
algebra by $ux+y$ where $u\in\bF_4$ is a primitive cubic root of unity.
\end{exam}

\subsection{Definition of $\fAlg_{n,\le d}$ and $\fAlg_{n,d}$}

Given a basis, a unitary algebra structure is determined by the
coefficients $c^{i,j}_k$ for all $2\le i,j\le n$ and $1\le k\le n$,
such that $e_ie_j=\sum_k\,c^{i,j}_ke_k$, subject to the relations
obtained by expansion of the associativity equations
$(e_ie_j)e_k=e_i(e_je_k)$. These relations tell us that
the moduli scheme $\Alg_n$ is a closed subscheme of affine space of
dimension $n(n-1)^2$ over $\bZ$. The change of basis is expressed by
an action of the subgroup $G\subset\GL_n$ which is the stabilizer of
the first basis vector, and $\fAlg_n$ is the quotient stack of
$\Alg_n$ by $G$.

We denote by $\fAlg_{n,\le d}$, resp. $\Alg_{n,\le d}$, the closed substack,
resp. the closed subscheme classifying $n$-dimensional algebras with degree of
algebraicity $\le d$.
We introduce also the locally open substacks
$$\fAlg_{n,d}=\fAlg_{n,\le d}\setminus \fAlg_{n,\le d-1}
\qquad\mbox{and}\qquad
\Alg_{n,d}=\Alg_{n,\le d}\setminus\Alg_{n,\le d-1}\;.$$
Our main goal in the next section is to construct determinant functions on the
normalized strata $\fAlg_{n,d}^\sim$ , by which we mean, the normalization of
the reduced stack $(\fAlg_{n,d})_\red$. We observe that $\fAlg_{n,d}^\sim$
is the quotient of $\Alg_{n,d}^\sim$ by $G$, because $G$ is smooth
and normality is local in the smooth topology.

\section{Determinants} \label{determinants}

\subsection{Construction of determinants}

Now comes our main result~:

\begin{theo} \label{det}
Let $S$ be a normal integral scheme, $A$ an $n$-dimensional
$\bG_{a,S}$-algebra scheme, and~$d\ge 2$ the least integer such that the
morphism $A\to\wedge^{d+1}A$ defined by $a\mapsto 1\wedge a\wedge\dots\wedge
a^d$ is zero. Then there exists a unique section of $\Sym^d(A^\vee)$, called the
{\em determinant of $A$} and denoted $\det_{A/S}$ or simply $\det$, such that~:
\begin{trivlist}
\itemn{1} $\det\colon A\to\bG_{a,S}$ is a morphism of multiplicative
unitary monoid schemes, and
\itemn{2} the morphism of schemes $A\to A$ defined by evaluating on
$a\in A$ the polynomial $P(T)=\det(T-a)$ is zero. (Here $T$ is a
scalar in the $\bG_{a,S}[T]$-algebra $A\times_S\bG_{a,S}[T]$.)
\end{trivlist}
The determinant satisfies the further properties~:
\begin{trivlist}
\itemn{3} the group scheme of units of $A$ is the
preimage of the multiplicative group $\bG_{m,S}$ under $\det$~;
\itemn{4} the formation of $\det$ is compatible with flat extensions of
normal integral schemes $S'\to S$.
\end{trivlist}
\end{theo}

Here are a few comments before we pass to the proof. One strong point of the
theorem is that $d$ is the degree of the generic fibre of $A/S$, not of all
fibres. If we want to extend it to more general base schemes, we can expect
such a strong result only for integral bases. An extension to arbitrary base
schemes may be reasonable if one assumes really $\deg(A/S)=d$ in the sense of
\ref{deg_alg} (i.e. all fibres of $A/S$ have degree $d$).
Intuitively, the basic idea is to produce the determinant as the constant
coefficient of the minimum polynomial of a generic element of the algebra,
just like in the case of a matrix algebra. More precisely, we consider the
universal element $\alpha$ of $A$, and we study the kernel of the evaluation
morphism $\ev_\alpha\colon \bG_{a,A}[T] \to A_A$. The left regular
representation provides an embedding of $A_A$ into $\End_A(A_A)$,
its algebra of module endomorphisms. Hence $\alpha$ is cancelled by its
Cayley-Hamilton polynomial $\CH_\alpha$. This provides a canonical polynomial
of degree $n$ in the kernel of $\ev_\alpha$. Over a normal integral scheme, we
prove that the kernel of $\ev_\alpha$ is generated by a single monic
polynomial. Over a general base scheme $S$, this is not to be expected~: for
instance, in example~\ref{exam1} two generators are required, namely
$\CH_\alpha$ and the polynomial $P(T)=\epsilon (T-r)^2$, where we have written
$\alpha=r+sx+ty$. So this raises the problem of finding a "distinguished"
polynomial in $\ker(\ev_\alpha)$. Ideally it should be
monic, it should divide $\CH_\alpha$ and have minimal degree.

We come to the proof of theorem~\ref{det}. Below, whenever we reduce to
the local case, we will use the following notation~: $S=\Spec(R)$ is affine
and small enough so that the algebra of sections $\sA=A(R)$ is free over $R$~;
the function ring of $A$ is $R_\sA=\Sym(\sA^\vee)$~;  the universal element is
$\alpha\in \sA\otimes_R R_\sA$. We denote by $K$ resp. $K_\sA$ the
fraction field of $R$ resp. $R_\sA$.

\begin{prooU}
Normality is not useful here~; it is enough to assume that $S$ is
integral. If~$S$ is irreducible, the integer $d$ defined in the statement of
the theorem is the same for all open subschemes of $S$,
hence proving unicity is a local question and we may assume that $S$ is
affine. We use the above notations for the affine case. Let $P\in K_\sA[T]$
be the minimum polynomial of the universal element $\alpha$. Clearly its
degree is~$d$. Consider the polynomial $Q(T)=\det_A(T-\alpha)$. Since the
determinant is required to take $1$ to $1$, and is homogeneous of
degree $d$, it follows that $Q$ is nonzero and monic. By property~(2) we have
$Q(\alpha)=0$, hence we get $Q=P$. Therefore $\det_A(\alpha)=(-1)^dP(0)$
and this determines $\det_A$.
\end{prooU}

We add that a flat extension of normal integral rings $R\to R'$ preserves the
kernels, hence the minimum polynomial of the universal element. Statement (4)
of the theorem follows easily, by immediate globalization. It remains to prove
that there exists a section $\det_{A/S}$ with the properties (1)-(2)-(3).

\begin{prooE}
By unicity, the question of existence is local and we may assume that
$S$ is affine and $\sA=A(R)$ is free over $R$. Since the minimum polynomial
$P$ divides the Cayley-Hamilton polynomial~$\CH_\alpha$, its roots (in some
algebraic closure of $K_\sA$) are roots of $\CH_\alpha$, hence integral
over~$R_\sA$. Therefore the coefficients of $P$ are themselves integral,
and since they lie in $K_\sA$ and $R_\sA$ is normal, it follows that
the coefficients of $P$ are in $R_\sA$.
We define $\det=\det_A:=(-1)^dP(0)$.
This is an element in the graded algebra $R_\sA$, and we now prove
that it lies in the degree $d$ component. In order to do so we
fix a basis $e_1=1,e_2,\dots,e_n$ for $\sA$. Let $t_i=e_i^*$
and write $\alpha=t_1e_1+\dots+t_ne_n$. Write the minimum polynomial of
$\alpha$ as $P(t_1,\dots,t_n,T)$, so that by definition
$$
P(t_1,\dots,t_n,t_1e_1+\dots+t_ne_n)=0\;.
$$
Let $X$ be an indeterminate. If we substitute $Xt_i$ to $t_i$ we find
$$
P(Xt_1,\dots,Xt_n,Xt_1e_1+\dots+Xt_ne_n)=0 \;,
$$
thus $P(Xt_1,\dots,Xt_n,XT)$ cancels $\alpha$. By unicity of the minimum
polynomial in $K_\sA(X)[T]$, it follows that
$P(Xt_1,\dots,Xt_n,XT)=X^dP(t_1,\dots,t_n,T)$. Therefore $P$ is homogeneous of
degree~$d$ in $t_1,\dots,t_n,T$ and hence $\det$ is homogeneous of
degree~$d$ in $t_1,\dots,t_n$.

We proceed to prove (1). For this we extend the base to the normal
domain $R':=R_{\sA\times \sA}=R_\sA\otimes_R R_\sA$, and we have to
prove that $\det(\alpha_1\alpha_2)=\det(\alpha_1)\det(\alpha_2)$ where
$(\alpha_1,\alpha_2)$ is the universal element of $\sA\times \sA$.
To prove this identity, we may embed $\sA\otimes_R R'$ into
$\sA\otimes_R \bar K{}'$, where~$\bar K{}'$ is an
algebraic closure of the fraction field of $R'$. For simplicity we
will write $K$ for $\bar K{}'$ and $\sA$ for $\sA\otimes_R \bar K{}'$.
Now let $V_0=0\subset V_1\subset V_2\subset\dots\subset V_d=A$ be a composition
series for the left $\sA$-module~$\sA$. Denote by $W_i$ the simple $\sA$-module
$V_{i}/V_{i-1}$. Its commutant $\End_\sA(W_i)$ is a division algebra~; since $K$
is algebraically closed, we have $\End_\sA(W_i)=K$. By Burnside's
theorem (\cite{B}, \S~4, n$^\circ$~3, corollaire~1), the morphism
$\varphi_i$ from $A$ to the bicommutant $\End_K(W_i)$ is surjective. It follows
from \ref{universal_elements} that the image in $\End_K(W_i)$ of the universal
element $\alpha$ is the universal element of $\End_K(W_i)$, and hence
its minimum polynomial $\chi_i$ is irreducible. If we choose a
$K$-basis of $\sA$ adapted to the composition series
$V_i$, then the regular left representation of $\sA$ provides an embedding
$\sA\hookrightarrow \End_K(\sA)$ whose image is block-triangular, with the $i$-th
diagonal block isomorphic to $\End_K(W_i)$. By the
Cayley-Hamilton theorem and the computation of a triangular determinant, the
image of $\alpha$ in $\End_K(\sA)$ is cancelled by the product
$\chi_1\dots\chi_n$. It follows that $P$ is a product of some of the
$\chi_i$'s, say the first $h$, so that
$\det=\det_1\dots\det_h$. The facts that $\det$ is multiplicative
and respects the unit follow.

We now prove (2). Let $Q(T)=\det(T-\alpha)$~; it remains to prove that
$Q(T)=P(T)$. For this it is enough to prove that $Q$ cancels $\alpha$. We
consider the algebra $\sA\otimes_R R'[T]$ and the element $T-\alpha$
therein.
If we write $P(T)=(-1)^dTP_0(T)+P(0)$ where $P_0$ has degree $d-1$, the equality
$P(\alpha)=0$ gives $\alpha P_0(\alpha)=\det(\alpha)$.
Since $\alpha$ is the universal element, we may look at this equality
in $\sA\otimes_R R'[T]$ and substitute $T-\alpha$ to $\alpha$. Thus~:
$$
(T-\alpha)\cdot P_0(T-\alpha)={\det}(T-\alpha)=Q(T)\;,
$$
Writing $P_0(T-\alpha)=T^{d-1}+\pi_1T^{d-2}+\dots+\pi_{d-1}$ and
$Q(T)=q_0T^d+\dots+q_d$ with $\pi_i,q_j\in \sA\otimes_R R'$, 
the above equality yields $\pi_0=q_0$, $\pi_1-\alpha \pi_0=q_1$,
..., $-\alpha \pi_{d-1}=q_d$. Now multiply the first equality by $1$, the
second by $\alpha$, ..., the last by $\alpha^d$ and add, we find $Q(\alpha)=0$.

It remains to prove (3), but this is obvious from the universal formula
$\alpha P_0(\alpha)=\det(\alpha)$.
\end{prooE}

\begin{prop} \label{det_on_Alg_d}
Let $\fA$ be the universal $n$-dimensional algebra of degree $d$ over
the normalized stratum $\fAlg_{n,d}^\sim$. Then $\fA$ has a determinant
satifying the properties (1) to (4) of theorem~\ref{det}, and its formation
commutes with any base change.
\end{prop}

\begin{proo}
The algebraic stack $\fAlg_{n,d}^\sim$ has an atlas
$\Alg_{n,d}^\sim\to\fAlg_{n,d}^\sim$ which is smooth and affine, hence the
fibre square of $\Alg_{n,d}^\sim$ over $\fAlg_{n,d}^\sim$ is a normal affine
scheme, denoted $\Alg^{(2)}$. By theorem~\ref{det} there exists a
determinant for $\fA$ after base extension to $\Alg_{n,d}^\sim$. The pullbacks
of this determinant via the two projections $\Alg^{(2)}\to \Alg_{n,d}^\sim$
are two functions satisfying all the properties of theorem~\ref{det}, hence by
unicity they coincide, and it follows that the determinant descends to
$\fAlg_{n,d}^\sim$.
\end{proo}

It follows from this proposition that for any $S$-algebra scheme $A$ such that
the corresponding classifying morphism $S\to \fAlg_n$ factors through
$f\colon S\to\fAlg_{n,d}^\sim$, we can define $\det_A:=f^*\det_\fA$.
When $S$ is the spectrum of a normal domain, this is the same as the
determinant given by  theorem~\ref{det}.

But of course, we would like more. The assumption of normality is used in one
single place, in the beginning of the existence part of the proof. Thus, it is
reasonable to expect that theorem~\ref{det} extends to an integral base scheme
$S$. Then we may pass to $(\fAlg_{n,d})_\red$ by working on the irreducible
components and then glue. So we ask~:

\begin{ques}
Can one construct a determinant on the reduced substack $(\fAlg_{n,d})_\red$~?
\end{ques}

\subsection{Properties of the determinant}

The determinant satisfies a strong invariance property with respect to
automorphisms~:

\begin{prop} \label{invariance_of_det}
Under the assumptions of theorem~\ref{det}, let $f\colon A\to A$ be a morphism
of $S$-schemes that is either a ring scheme automorphism or a ring scheme
antiautomorphism, and takes the scalars to scalars
(i.e. $f(\bG_{a,S})=\bG_{a,S}$). Then $\det_{A/S}\circ f=f\circ\det_{A/S}$. In
particular if $f$ is a $\bG_{a,S}$-algebra (anti)automorphism then
$\det_{A/S}\circ f=\det_{A/S}$.
\end{prop}

\begin{proo}
Let $\alpha$ be the universal element and let $P(T)=\det_{A/S}(T-\alpha)$ be its
minimum polynomial. If we let $f$ act on the polynomials by $f(T)=T$
and its natural action on the coefficients, since it is a ring
(anti)automorphism we find that $fP$ is the minimum polynomial of
$f(\alpha)$. Moreover, $f(\alpha)$ is also a
universal element, so it is uniquely a pullback of $\alpha$, indeed, via
$f$. Accordingly $fP$ is the pullback $f^*P$. Looking at the constant term, we
get the result.
\end{proo}

Here is another useful property, which for simplicity we state in the
affine case~:

\begin{prop} \label{det_irred}
Let $R$ be a normal domain and $\sA$ an $n$-dimensional $R$-algebra.
If the determinant $\det_{\sA/R}$ is irreducible, then its degree $d$
divides $n$.
\end{prop}

\begin{proo}
Let $\alpha$ be the universal element and let
$P(T)=\det_{\sA/S}(T-\alpha)$ be its minimum polynomial. If we extend
the natural degree of $R_\sA=\Sym(\sA^\vee)$ by assigning degree $1$
to $T$, then $P$ is homogeneous of degree $d$. So if $P=QR$ in
$R_\sA[T]$, then $Q$ and $R$ are homogeneous of degrees $e$, $f$ such
that $d=e+f$. Since we assumed that $\det_{\sA/R}=(-1)^dQ(0)R(0)$ is
irreducible, it follows that either $Q(0)\in A^\times$ and $e=0$, or
$R(0)\in A^\times$ and $f=0$. Therefore $P$ is irreducible in
$R_\sA[T]$. Moreover the Cayley-Hamilton polynomial of $\alpha$ has
the same irreducible factors as $P$, so by irreducibility
$\CH_\alpha=P^m$ for some $m$. By taking degrees we get $n=dm$.
\end{proo}

\section{Computations and applications} \label{comp_appl}

\subsection{More determinants}

We start with the determinant of the opposite algebra, and the determinant of
a product.

\begin{prop} \label{det_of_A_opp}
Under the assumptions of theorem~\ref{det}, let $A^\circ\to S$ be
the opposite algebra whose multiplication is the opposite as that of $A$,
i.e. $a\star_{A^\circ} b=b\star_A a$. Then $\det_{A^\circ/S}=\det_{A/S}$.
\end{prop}

\begin{proo}
We have $A^\circ=A$ as vector bundles, so
$\Sym^d(A^{\circ\vee})=\Sym^d(A^\vee)$. It is clear that $\det_{A/S}$
satisfies the properties of theorem~\ref{det} for $A^\circ$, hence by unicity
$\det_{A^\circ/S}=\det_{A/S}$.
\end{proo}

\begin{prop} \label{product}
Let $S$ be a normal integral scheme.
\begin{trivlist}
\itemn{1} If $f\colon A\to B$ a surjective homomorphism of finite-dimensional
$\bG_{a,S}$-algebra schemes with determinants $\det_{A/S}\in\Sym^d(A^\vee)$ and
$\det_{B/S}\in\Sym^e(B^\vee)$, then there exists a unique section $\det_{A/B/S}$
of $\Sym^{d-e}(A^\vee)$ such that $\det_{A/S}=\det_{A/B/S}\cdot
f^*\det_{B/S}$. Moreover $\det_{A/B/S}\colon A\to\bG_{a,S}$ is a
morphism of multiplicative unitary monoid schemes.
\itemn{2} If $C=A\times B$ is a product algebra, then
$\det_{C/S}=\pr_1^*\det_{A/S}\cdot\pr_2^*\det_{B/S}$.
\end{trivlist}
\end{prop}

\begin{proo}
(1) We can work locally over the base and hence suppose that $S=\Spec(R)$ is
affine. Consider the $R$-algebra $\sA=A(R)$, the function ring
$R_\sA=\Sym(\sA^\vee)$, its fraction field $K_\sA=\Frac(R_\sA)$, and the
universal element $\alpha$. Similarly we have $\sB$, $R_\sB$, $K_\sB$,
$\beta$.

Let $P\in R_\sA[T]$ and $Q\in R_\sB[T]$ be the minimum polynomials of $\alpha$
and $\beta$. The relation $(f\times\Id_A)\circ\alpha=f^*\beta$ (see
\ref{image_of_univ_elt}) applied to $P(\alpha)=0$ gives $P(\beta)=0$. Since
$f$ is a surjective map of vector bundles, it is dominant~; in particular
$f^*\colon R_\sB\to R_\sA$ is injective. It follows that $Q$ divides $P$ in
$K_\sA[T]$, that is to say $P=QR$ for some
$R\in K_\sA[T]$. In fact $R\in R_\sA[T]$ since $Q$ is unitary. We define
$\det_{A/B/S}:=(-1)^{d-e}R(0)$. It is clear that it satisfies all the properties
asserted in the statement of the proposition.

\no (2) If $C=A\times B$, we have projections
$\pr_1\colon C\to A$ and $\pr_2\colon C\to B$. We localize and suppose that
$S=\Spec(R)$ is small enough so that $A$ and $B$ are trivial as vector
bundles. We end up with $\sC=\sA\times\sB$, $R_\sC=R_\sA\otimes R_\sB$,
and the fraction field $K_\sC$. With a slight abuse of notation we can write
the universal element of $C$ as $\gamma=(\alpha,\beta)$, if we think of the
injections $\pr_1^*$  and $\pr_2^*$ as inclusions. From the above,
the minimum polynomial of $\gamma$ is a multiple of $P$ and $Q$.
We now argue that $P$ and $Q$ are coprime in $K_\sC[T]$. For otherwise, the
resultant $\Res(P,Q)$ vanishes, giving a relation of algebraic dependance
between the variables of $R_\sA$ and $R_\sB$. Since $R_\sA$ and~$R_\sB$ are
polynomial rings in independent variables, this relation can
involve only elements from~$R$, hence $P$ and $Q$ belong to
$R[T]$. But this is impossible since $d\ge 1$. Therefore $P$ and $Q$ are
coprime, so the minimum polynomial of $\gamma$ is $PQ$, and
$\det_{C/S}=\det_{A/S}\det_{B/S}$.
\end{proo}

It is natural to ask what is the dual result, namely, what is the determinant
of a coproduct algebra. However the coproduct in the category of algebras is
the free product, which is {\em not} finite-dimensional. Still, one can wish
to compute the determinant of a tensor product $A\otimes B$, which is
universal for pairs of maps $f:A\to C$, $g:B\to C$ whose images $f(A)$ and
$g(B)$ commute. If $f\colon A\to B$ is an injective homomorphism of algebra
schemes over a normal base, then it is not hard to prove that $\det_A$ divides
$\det_B|_A$. In the case of a tensor product we get that $\det_A$ divides
$\det_{A\otimes B}|_A$ and $\det_B$ divides $\det_{A\otimes B}|_B$, but it is
not obvious how to guess an expression for $\det_{A\otimes B}$. For
example, in view of the isomorphism $\MM_p(k)\otimes_k\MM_q(k)\simeq
\MM_{pq}(k)$, the formula we are
looking for should give the determinant of $(pq,pq)$ matrices in terms
of the determinant of $(p,p)$ and $(q,q)$ matrices. To sum up~:

\begin{ques}
Is there a simple ``formula'' for the determinant of a tensor product
$A\otimes B$ in terms of $\det_A$ and $\det_B$~?
\end{ques}

\subsection{Traces and discriminants} \label{tr_and_disc}

Let $S$ be a normal integral scheme, $A$ an $n$-dimensional
$\bG_{a,S}$-algebra scheme, $\det_A$ its determinant of degree $d$, and
$\alpha$ the universal element. We define some invariants of the intrinsic
structure of $A$.

\begin{noth}
{\em Coefficients of the characteristic polynomial}. They are
the sections $c_i$ of $\Sym^i(A^\vee)$ defined by
$\det_A(T-\alpha)=T^d-c_1T^{d-1}+\dots+(-1)^dc_d$.
The coefficient $c_1$ is called the {\em trace} and denoted $\tr_A$.
From theorem~\ref{det} it follows that the formation of the $c_i$ commutes
with flat extensions of normal integral schemes.
\end{noth}

\begin{noth}
{\em Discriminant}.
Locally over $S$ we may choose a basis $e_1,\dots,e_n$ for $A$ and compute the
determinant of the matrix whose $(i,j)$-th element is $\tr_A(e_ie_j)$. By the
usual formula, a base change multiplies this by the square of the determinant
of the transition matrix. Therefore we obtain a section of the quotient
(monoid) stack $[\bG_{a,S}/\bG_{m,S}]$, where the multiplicative group
$\bG_{m,S}$ acts on $\bG_{a,S}$ by the rule $z.x=z^2x$. These local sections
are canonical and hence glue to a section over all of $S$. This is called the
{\em discriminant} of $A$ and denoted $\disc(A)$. Here again, the formation of
$\disc(A)$ commutes with flat extensions of normal integral schemes.
\end{noth}

\begin{noth}
{\em Unimodular group}. The {\em unimodular group} of $A/S$ is the kernel of
the determinant, namely $U=(\det_A)^{-1}(1)$. This is a subgroup scheme of the group
$G=A^\times=(\det_A)^{-1}(\bG_{m,S})$ of invertible elements of $A$.
\end{noth}

\subsection{Application to the topology of $\fAlg_n$} \label{appl_top}

In this subsection, the base ring is a field $k$, and the moduli stack
$\fAlg_n$ is considered over $k$.

\begin{prop}
$\Alg_{n,n}$ is irreducible of dimension $n^2$, so $\overline{\Alg}_{n,n}$ is
an irreducible component of $\Alg_n$. In other words, $\Alg_{n,\le n-1}$ is
a union of irreducible components of $\Alg_n$.
\end{prop}

\begin{proo}
By lemma~\ref{monogenic}, $\fAlg_{n,n}$ is the open substack of algebras that
are locally (for the {\'e}tale topology) monogenic, hence
commutative. Therefore the result follows from \cite{P}, \S~6. The
argument is so simple that we recall it shortly~: the locus of {\'e}tale algebras
$E\subset \fAlg_{n,n}$ is defined by the nonvanishing of the discriminant,
hence it is open and dense. Since an {\'e}tale algebra splits after an {\'e}tale
base extension, the orbit of the split algebra $A=k^n$ under $\GL_n$ acting by
base change is $E$. It follows that $E$ is irreducible, hence so is
$\fAlg_{n,n}$.
\end{proo}

In \cite{P}, questions 9.8-9.9, B. Poonen asks what is the functor of points
of $\overline{\Alg}_{n,n}$ (his notation for $\overline{\Alg}_{n,n}$ is
$\overline{\fB_n^\et}$). We have no complete answer, but our result indicates
that all monogenic algebras are points of this functor.

Now here is an application of the determinant to the topology of $\fAlg_n$~:

\begin{prop} \label{an_application}
Let $n_1,\dots,n_r\ge 1$ be integers and $n=(n_1)^2+\dots+(n_r)^2$. Assume
that one of the $n_i$ is at least $2$ and denote by $\nu$ their infimum.
Then the irreducible component of the algebra
$\sA_0=\MM_{n_1}(k)\times\dots\times\MM_{n_r}(k)$ is contained in
$\fAlg_{n,\le n-\nu}$.
\end{prop}

\begin{proo}
Assume that the irreducible component of $\sA_0$ meets $\fAlg_{n,d}$.
Then there is a discrete valuation ring $R$ with fraction field $K$ and
residue field $k$, and an $R$-algebra $\sA$ such that $\sA\otimes k\simeq
\sA_0$ and $\sA\otimes K$ has degree $d$. Let $\alpha$ (resp. $\alpha_0$)
denote the universal element of $\sA$ (resp. $\sA_0$). If we had $d=n$, then
$\sA\otimes K$ would be commutative by lemma~\ref{monogenic}, and then $\sA_0$
would be commutative also. Since this is not the case, we have $d<n$. We have
a factorization of the Cayley-Hamilton polynomial of $\alpha$ as
$\CH_\alpha(T)=\psi(T)\det_\sA(T-\alpha)$, where the degree of $\psi$ is
$n-d>0$. There is an analogous factorization for the Cayley-Hamilton
polynomial of $\alpha_0$, which is just the reduction of $\CH_\alpha(T)$
modulo the maximal ideal of $R$. Let $\delta_i$ be the determinant of the
algebra $\MM_{n_i}(k)$, irreducible of degree $n_i$. We know from
proposition~\ref{product} that $\det_{\sA_0}(T-\alpha_0)$ is the product of
the $\delta_i$. Moreover $\CH_{\alpha_0}(T)$ and $\det_{\sA_0}(T-\alpha_0)$
have the same irreducible factors, hence $\bar\psi(T)$ is divisible by one of
the $\delta_i$. Therefore its degree is at least $\nu$, that is, $n-d\ge
\nu$. So $d\le n-\nu$.
\end{proo}

These results raise the natural question~:

\begin{ques} \label{ques_irred_comp}
Is $\fAlg_{n,\le d}$ (or $\Alg_{n,\le d}$) a union of irreducible
components of $\fAlg_n$ (or $\Alg_n$)~?
\end{ques}

\section{Examples} \label{examples}

\subsection{Determinants of some classical algebras}

\begin{trivlist}
\itemn{1} Let $K/k$ be a finite Galois field
extension with Galois group $G$. Then $\det_{K/k}$ is the norm defined by
$N(x)=\prod_{g\in G}\,g(x)$.
\itemn{2} Let $k=k_0(t_1,\dots,t_s)$ be the field of rational functions in
$t_1,\dots,t_s$ over a field of characteristic $p>0$ and let $K/k$ be the
finite purely inseparable extension generated by $s$ elements
$u_i=(t_i)^{1/p}$. A natural basis of $K/k$ is given by the monomials
$u^I:=(u_1)^{i_1}\dots (u_s)^{i_s}$ for $0\le i_1,\dots,i_s\le p-1$.
For $x=\sum_I\,x_Iu^I$ one finds $\det_{K/k}(x)=-\sum_I\,(x_I)^pt^I$.
(Note that $x^p=\sum_I\,(x_I)^pt^I$).
\itemn{3} Let $H=H_{\alpha,\beta}$ be the quaternion algebra generated by
$1,i,j,k$ with relations $i^2=\alpha$, $j^2=\beta$, $ij=-ji=k$. Then $\det_{H/k}$ is the
norm $N(a+bi+cj+dk)=a^2-\alpha b^2-\beta c^2+\alpha\beta d^2$.
\itemn{4} Let $A=\wedge^* E$ be the exterior algebra of an
$r$-dimensional vector space $E$. To any basis $\{e_1,\dots,e_n\}$ of
$E$ is associated a basis of $A$ whose $I$-th vector is
$e_I=e_{i_1}\wedge\dots\wedge e_{i_k}$, for all $0\le k\le n$ and
$I=\{i_1<\dots<i_k\}\subset \{1,\dots,n\}$. Then
$\det_{A/k}(x)=(x_\emptyset)^r$, where $x=\sum\, x_Ie_I$.
\end{trivlist}

It would be interesting to generalize examples (3) and (4) by computing the
determinant of a general Clifford algebra.


%

\subsection{Algebras of dimension $3$}

The simplest case is dimension $n=2$, but then the determinant is just the
Cayley-Hamilton polynomial. Precisely, for a $2$-dimensional algebra endowed
with a basis $\{1,x\}$, the datum of an equation $x^2=ax+b$
determines a multiplication that is automatically associative, thus $\Alg_2$
is an affine plane. One checks that $\det(r+sx)=r^2-s^2b+ars$.

Next there is the case of dimension $n=3$, which is still easily
tractable. In~\cite{MT}, Miranda and Teicher study the noncommutative
$3$-dimensional algebras over integral schemes. It is in fact not more
complicated to decribe directly the whole moduli stack
$\fAlg_3$. As we will see, the scheme $\Alg_3$ has two irreducible
components~: one is the component of commutative algebras, which is
also the closure of $\Alg_{3,3}$, and the other is just $\Alg_{3,\le
  2}=\Alg_{3,2}$. We will pay particular attention to
$\Alg_{3,2}$, by describing the base change group action and giving
the determinant on $\Alg_{3,2}$.

\begin{noth} {\em Equations for $\Alg_3$.}
The multiplication table of a $3$-dimensional algebra with basis
$\{1,x,y\}$ looks like this~:
$$
\begin{array}{ccl}
x^2 & = & a+bx+cy \\
xy & = & d+ex+fy \\
yx & = & g+hx+iy \\
y^2 & = & j+kx+ly \\
\end{array}
$$
The conditions of associativity are given by the following eight sets of
equations~:

\bigskip

\begin{center}
\begin{tabular}{ccc}
\begin{tabular}{|c|}
\hline
\begin{minipage}{4.5cm}
$$(x^2)x=x(x^2)$$
\end{minipage} \\
\hline
$c(d-g)=0$ \\
$c(e-h)=0$ \\
$c(f-i)=0$ \\
\hline
\end{tabular}
& &
\begin{tabular}{|c|}
\hline
\begin{minipage}{4.5cm}
$$(x^2)y=x(xy)$$
\end{minipage} \\
\hline
$bd+cj=ea+fd$ \\
$be+ck=d+eb+fe$ \\
$a+bf+cl=ec+f^2$ \\
\hline
\end{tabular} \\
\end{tabular}
\end{center}

\begin{center}
\begin{tabular}{ccc}
\begin{tabular}{|c|}
\hline
\begin{minipage}{4.5cm}
$$(xy)x=x(yx)$$
\end{minipage} \\
\hline
$ea+fg=ha+id$ \\
$d+eb+fh=g+hb+ie$ \\
$ec+fi=hc+if$ \\
\hline
\end{tabular}
& &
\begin{tabular}{|c|}
\hline
\begin{minipage}{4.5cm}
$$(xy)y=x(y^2)$$
\end{minipage} \\
\hline
$ed+fj=ka+ld$ \\
$e^2+fk=j+kb+le$ \\
$d+ef+lf=kc+lf$ \\
\hline
\end{tabular} \\
\end{tabular}
\end{center}

\begin{center}
\begin{tabular}{ccc}
\begin{tabular}{|c|}
\hline
\begin{minipage}{4.5cm}
$$(yx)x=y(x^2)$$
\end{minipage} \\
\hline
$ha+ig=bg+cj$ \\
$g+hb+ih=bh+ck$ \\
$hc+i^2=a+bi+cl$ \\
\hline
\end{tabular}
& &
\begin{tabular}{|c|}
\hline
\begin{minipage}{4.5cm}
$$(yx)y=y(xy)$$
\end{minipage} \\
\hline
$hd+ij=eg+fj$ \\
$he+ik=eh+fk$ \\
$g+hf+il=d+ei+fl$ \\
\hline
\end{tabular} \\
\end{tabular}
\end{center}

\begin{center}
\begin{tabular}{ccc}
\begin{tabular}{|c|}
\hline
\begin{minipage}{4.5cm}
$$(y^2)x=y(yx)$$
\end{minipage} \\
\hline
$ka+lg=hg+ij$ \\
$j+kb+lh=h^2+ik$ \\
$kc+li=g+hi+il$ \\
\hline
\end{tabular}
& &
\begin{tabular}{|c|}
\hline
\begin{minipage}{4.5cm}
$$(y^2)y=y(y^2)$$
\end{minipage} \\
\hline
$k(d-g)=0$ \\
$k(e-h)=0$ \\
$k(f-i)=0$ \\
\hline
\end{tabular} \\
\end{tabular}
\end{center}

\bigskip

These twenty-four equations can be rearranged and simplified to give the
following twelve~:
$$
\left\{
\begin{array}{l}
a=f(f-b)+c(e-l) \\
d=ck-ef \\
g=ck-hi \\
j=k(f-b)+e(e-l) \\
\end{array}
\right.
\qquad\mbox{and}\qquad
\left\{
\begin{array}{l}
c(e-h)=c(f-i)=0 \\
k(e-h)=k(f-i)=0 \\
(f+i-b)(e-h)=(f+i-b)(f-i)=0 \\
(e+h-l)(e-h)=(e+h-l)(f-i)=0
\end{array}
\right.
$$
The first four equations imply that $\Alg_3$ is a subscheme of affine
$8$-space in coordinates $b$, $c$, $e$, $f$, $h$, $i$, $k$, $l$. The remaining
eight equations generate an ideal which is the product of the two prime ideals
$p_1=(e-h,f-i)$ and $p_2=(c,k,f+i-b,e+h-l)$. So $\Alg_3$ is defined in $\bA^8$
by the ideal $p_1p_2$, it has two irreducible components which are affine
spaces.
\end{noth}

\begin{noth} {\em The component of commutative algebras.}
The component defined by $p_1=0$ is the component of commutative algebras~; it
is also the closure of $\Alg_{3,3}$~; it is isomorphic to affine $6$-space
with coordinates $b,c,e,f,k,l$. The multiplication table of an algebra in this
component is the following~:
$$
\begin{array}{ccl}
x^2 & = & f(f-b)+c(e-l)+bx+cy \\
xy = yx & = & ck-ef+ex+fy \\
y^2 & = & k(f-b)+e(e-l)+kx+ly \\
\end{array}
$$
\end{noth}

\begin{noth} {\em The component of degree $2$ algebras.}
We now pay more attention to the component defined by $p_2=0$, and we start by
checking that it is $\Alg_{3,\le 2}=\Alg_{3,2}$. Let $\alpha=r+sx+ty$ be
the universal element, we have
$$
\alpha^2=\big[r^2+st(d+g)+s^2a+t^2j\big]
+\big[2rs+st(e+h)+s^2b+t^2k\big]x+\big[2rt+st(f+i)+s^2c+t^2l\big]y\;.
$$
The locus $\Alg_{3,\le 2}$ is defined by the vanishing of
$1\wedge\alpha\wedge\alpha^2$, hence~:
$$
\left|
\begin{array}{cc}
2rs+st(e+h)+s^2b+t^2k & s \\
2rt+st(f+i)+s^2c+t^2l & t
\end{array}
\right|\;=-cs^3+(b-f-i)s^2t+(e+h-l)st^2+kt^3=0
$$
as polynomials in $s,t$. We find $b=f+i$, $c=k=0$, $l=e+h$, in other
words this is indeed the component defined by $p_2=0$. Thus $\Alg_{3,2}$ is
isomorphic to affine $4$-space in coordinates $e,f,h,i$. The multiplication
table of an algebra in this component is the following~:
$$
\begin{array}{ccl}
x^2 & = & -fi+(f+i)x \\
xy & = & -ef+ex+fy \\
yx & = & -hi+hx+iy \\
y^2 & = & -eh+(e+h)y \\
\end{array}
$$
The minimum polynomial of $\alpha$ is $P(T)=(T-(r+is+et))\,(T-(r+fs+ht))$ and
the determinant is $\det(\alpha)=(r+is+et)(r+fs+ht)$. Note that by
proposition~\ref{det_irred}, we expected it to split.

We shortly describe $\fAlg_{3,2}$. Base change $x'=r+sx+ty$,
$y'=u+vx+wy$ is encoded by a transition matrix
$$
\left(
\begin{array}{ccc}
1 & r & u \\ 0 & s & v \\ 0 & t & w
\end{array}
\right)
$$
Set $\delta=sw-tv$ so that
$x=\frac{tu-rw}{\delta}+\frac{w}{\delta}x'-\frac{t}{\delta}y'$ and
$y=\frac{vr-su}{\delta}-\frac{v}{\delta}x'+\frac{s}{\delta}y'$. Then
one can calculate $(x')^2$, $x'y'$, $y'x'$, $(y')^2$ in terms of
$x'$, $y'$ and after a (tedious) computation one finds~:
$$
\begin{array}{l}
e'=u+vi+we \\
f'=r+sf+th \\
h'=u+vf+wh \\
i'=r+si+te \;.
\end{array}
$$
The stabilizer of a point of $\Alg_{3,2}$ with coordinates
$(e,f,h,i)$ has equations
$$
\begin{array}{l}
u+vi+(w-1)e=0 \\
r+(s-1)f+th=0 \\
u+vf+(w-1)h=0 \\
r+(s-1)i+te=0 \;,
\end{array}
$$
or in other words $v(f-i)=(w-1)(e-h)$, $(s-1)(f-i)=t(e-h)$,
$u+vi+(w-1)e=0$, $r+(s-1)f+th=0$. Let $F$ be the $2$-dimensional
intersection of $\Alg_{3,2}$ with the other irreducible component and
$U=\Alg_{3,2}\setminus F$ be the complement. We see that for points in
$F$, the stabilizer has dimension $4$. The image of $F$ in
$\fAlg_{3,2}$ is a single point which is the unique commutative
algebra, isomorphic to $k[x,y]/(x^2,xy,y^2)$, with automorphism group
$\GL_2(k)$. Away from $F$, the stabilizer has dimension $2$ and the
image of $U$ in $\fAlg_{3,2}$ has dimension $2$, composed of
algebras with $2$-dimensional automorphism group.
\end{noth}

\bigskip

\begin{flushleft}
{\em
Institut de Math{\'e}matiques

Th{\'e}orie des Nombres

Universit{\'e} Pierre et Marie Curie

Case 82

4, place Jussieu

F-75252 Paris Cedex 05

romagny@math.jussieu.fr}
\end{flushleft}


\begin{thebibliography}{99}
\bibitem[B]{B} {\sc N. Bourbaki}, {\it {\'E}l{\'e}ments de
    math{\'e}matique, Alg{\`e}bre, Chapitre 8~: Modules et anneaux
    semi-simples}, Actualit{\'e}s Sci. Ind. no. 1261, Hermann (1958).
\bibitem[DS]{DS} {\sc T. Dana-Picard, M. Schaps}, {\it Non-reduced
    components of $\Alg_n$}, Algebras and mo\-dules~II
  (Geiranger, 1996), 111--120, CMS Conf. Proc., 24,
  Amer. Math. Soc. (1998).
\bibitem[EGA]{EGA} {\sc A. Grothendieck}, {\it {\'E}l{\'e}ments de g{\'e}om{\'e}trie
    alg{\'e}brique II}, Publ. Math. Inst. Hautes {\'E}tudes Sci. No. 8 (1961).
\bibitem[Ga]{Ga} {\sc P. Gabriel}, {\it Finite representation type is
    open}, Proc. of the International Conference on Representations of
  Algebras (Carleton Univ., Ottawa, Ont., 1974), Carleton
  Math. Lecture Notes, No. 9 (1974).
\bibitem[LBR]{LBR} {\sc L. Le Bruyn, Z. Reichstein}, {\it Smoothness in
  algebraic geography}, Proc. London Math. Soc. 79 (1999) 158--190.
\bibitem[MT]{MT} {\sc R. Miranda, M. Teicher}, {\it Non-commutative algebras of
  dimension three over integral schemes}, Trans. AMS, v. 292, no. 2, (1985), 705--712.
\bibitem[P]{P} {\sc B. Poonen}, {\it The moduli space of commutative algebras
  of finite rank}, to appear in J. Europ. Math. Soc.
\end{thebibliography}
\end{document}